# A NOTE ON LASCAR STRONG TYPES IN SIMPLE THEORIES

BYUNGHAN KIM

July 25, 1996

ABSTRACT. Let $T$ be a countable, small simple theory. In this paper, we prove for such $T$, the notion of Lascar strong type coincides with the notion of strong type, over an arbitrary set.

§**1. Introduction.** We recall from [S1] and [S2] that a type $p$ *divides* over a set $A$ if there are a formula $\varphi(\bar{x}, \bar{y}) \in L$ and an $A$-indiscernible sequence $\langle \bar{a}_i | i \in \omega \rangle$ such that $p \vdash \varphi(\bar{x}, \bar{a}_0)$, and $\{\varphi(\bar{x}, \bar{a}_i) | i \in \omega\}$ is inconsistent. $p$ *forks* over $A$ if there are formulas $\varphi_0(\bar{x}, \bar{a}_0), ..., \varphi_n(\bar{x}, \bar{a}_n)$ such that $p \vdash \bigvee_{0 \leq i \leq n} \varphi_i(\bar{x}, \bar{a}_i)$, and $\varphi_i(\bar{x}, \bar{a}_i)$ divides over $A$ for each $i$. A first order theory $T$ is said to be *simple* if, for any complete type $p \in S(B)$, $p$ does not fork over some set $A(\subseteq B)$ with $|A| \leq |T|$. We also recall that $T$ is unstable if there are a formula $\varphi(\bar{x}, \bar{y})$ and tuples $\bar{a}_i$ ($i \in \omega$) such that $\models \varphi(\bar{a}_i, \bar{a}_j)$ if and only if $i \leq j \in \omega$. $T$ is said to be *stable* if $T$ is not unstable. It is well known that every stable theory is simple.

In [K1], it was proved that for simple $T$, the notion of forking is equivalent to that of dividing, and nonforking satisfies the following axioms.

(i) (*Symmetry*) $tp(\bar{a}/A\bar{b})$ does not fork over $A$ if and only if $tp(\bar{b}/A\bar{a})$ does not fork over $A$.

(ii) (*Transitivity*) Let $A \subseteq B \subseteq C$. Then $tp(\bar{a}/C)$ does not fork over $A$ if and only if $tp(\bar{a}/B)$ does not fork over $A$ and $tp(\bar{a}/C)$ does not fork over $B$.

(iii) (*Extension*) ([S2]) Let $p$ be a complete type in $S(A)$. For any set $B(\supseteq A)$, $p$ has a nonforking extension $q$ in $S(B)$.

(iv) (*Local Character*) For any complete type $p$ over $B$, there is a subset $A$ of $B$ such that $|A| \leq |T|$ and $p$ does not fork over $A$.

(v) (*Finite Character*) Let $A \subseteq B$. Then $tp(\bar{a}/B)$ does not fork over $A$ if and only if for each finite tuple $\bar{b} \in B$, $tp(\bar{a}/A\bar{b})$ does not fork over $A$.

In [KP], it was also shown that for simple $T$, nonforking satisfies the following additional axiom.

(vi) (*The Independence Theorem over a model*) Suppose that, for some model $M$, $tp(\bar{a}/M\bar{b})$ does not fork over $M$, and $p \in S(M)$. Let $p_1 \in S(M\bar{a})$ and $p_2 \in S(M\bar{b})$ be nonforking extensions of $p$. Then $p$ has a nonforking extension $p_3 \in S(M\bar{a}\bar{b})$ such that $p_1 \cup p_2 \subseteq p_3$.

Moreover, it was shown conversely that any theory equipped with an abstract relation between complete types and sets satisfying all the axioms (i) to (vi), must be simple, and the relation must be nonforking. This characterization of simple theories is a generalization of well known analogous fact on stable theories with axioms (i) to (v) and the uniqueness axiom.

(vi)$'$ (*Uniqueness*) Let $A \subseteq B$, and $p \in S(A)$. If either $A$ is a model, or an algebraically closed set in $\mathcal{C}^{eq}$, then $p$ has a unique nonforking extension in $S(B)$.

For stable $T$, the Independence Theorem over an algebraically closed set ( in $\mathcal{C}^{eq}$) obviously follows from the uniqueness axiom. (But the Independence Theorem does not hold over an arbitrary set. An equivalence relation having finitely many infinite equivalence classes trivially supplies a counterexample.) Therefore it is natural to ask for simple $T$, whether the Independence Theorem holds over an algebraically closed set in $\mathcal{C}^{eq}$. It is equivalent to ask whether the Independence Theorem holds for strong types. This issue will be discussed in this paper.

In fact, in [KP], the Independence Theorem for Lascar strong types was proved (for simple $T$), instead of strong types (Theorem 4). Here we recall the definition of Lascar strong type.

**Definition 1.** (i) By $Autf_A(\mathcal{C})$ we mean the subgroup of $Aut(\mathcal{C})$ generated by $\{f \in Aut(\mathcal{C})| f \in Aut_M(\mathcal{C}) \text{ for some model } M \supseteq A\}$.

(ii) Let $\bar{a}, \bar{b}$ be tuples of the same length. We say that $Lstp(\bar{a}/A) = Lstp(\bar{b}/A)$ ($\bar{a}$ and $\bar{b}$ have the same Lascar strong type over $A$) if there is $f \in Autf_A(\mathcal{C})$ such that $f(\bar{a}) = \bar{b}$.

We can think of the Lascar strong type over $A$ as simply the specification of an orbit under $Autf_A(\mathcal{C})$. Actually, $Lstp(\bar{a}/A) = Lstp(\bar{b}/A)$ if and only if there are tuples $\bar{a} = \bar{a}_0, \bar{a}_1, ..., \bar{a}_n = \bar{b}$ and models $M_1, M_2, .., M_n$, each containing $A$, such that $tp(\bar{a}_{i-1}/M_i) = tp(\bar{a}_i/M_i)$ for each $i = 1, .., n$. It follows that $Lstp(\bar{a}/A) = Lstp(\bar{b}/A)$ implies $stp(\bar{a}/A) = stp(\bar{b}/A)$, for any $T$. Moreover, for stable $T$, the Lascar strong type of $\bar{a}$ over $A$ is the same as the strong type of $\bar{a}$ over $A$, because $\bar{a}, \bar{b}$ realize the same strong type over $A$ if and only if $\bar{a}, \bar{b}$ realize the same type over a some model $M$ containing $A$. But, two notions are different in general. (See Example 12.)

Now let us state some preliminary facts from [KP].

**Proposition 2.** *Assume that $T$ is simple. For any $A \subseteq B$ and $\bar{a}$, there is $\bar{b}$ such that $Lstp(\bar{a}/A) = Lstp(\bar{b}/A)$ and $tp(\bar{b}/B)$ does not fork over $A$.*

**Proposition 3.** *Let $T$ be simple. Suppose that $Lstp(\bar{a}/A) = Lstp(\bar{b}/A)$ and $tp(\bar{a}/A\bar{b})$ does not fork over $A$. Then there is $M \supseteq A$ such that $tp(\bar{a}\bar{b}/M)$ does not fork over $A$, and $tp(\bar{a}/M) = tp(\bar{b}/M)$.*

**Theorem 4.** *Assume $T$ to be simple. Suppose $A, B, C, \bar{d}, \bar{e}$ satisfy*
  (i) $A \subseteq B$, $A \subseteq C$, $\{B, C\}$ *is $A$-independent,*
  (ii) $tp(\bar{d}/B)$ *does not fork over $A$, and $tp(\bar{e}/C)$ does not fork over $A$,*
  (iii) $Lstp(\bar{d}/A) = Lstp(\bar{e}/A)$.
*Then there is $\bar{a}$ such that $tp(\bar{a}/B \cup C)$ extends both $tp(\bar{d}/B)$ and $tp(\bar{e}/C)$, $tp(\bar{a}/B \cup C)$ does not fork over $A$, and $Lstp(\bar{a}/A) = Lstp(\bar{d}/A)(= Lstp(\bar{e}/A))$.*

The following proposition (in [K2]) is an easy consequence of the Independence Theorem over a model.

**Proposition 5.** *Let $T$ be simple. Let $\{\bar{a}_i | i \in \kappa\}$ be $M$-independent and $p \in S(M)$. Suppose that $p_i \in S(M\bar{a}_i)$ is a nonforking extension of $p$ for each $i \in \kappa$, then $\cup \{p_i | i \in \kappa\}$ is consistent and does not fork over $M$.*

In this paper, we summarize basic facts on Lascar strong types, and further investigate Lascar strong types and the Independence Theorem for simple theories. In particular we show that, for a countable small simple theory $T$, the notion of Lascar strong type is equivalent to that of strong type. This implies, in such $T$, the Independence Theorem holds over any algebraically closed sets in $\mathcal{C}^{eq}$.

The notation follows usual conventions. $T$ is a complete theory with no finite models in a first order language $L$. $p, q$ denote types, possibly partial. We work in a huge $\bar{\kappa}$-saturated model $\mathcal{C}$ as usual. Sets $A, B, C, ..$ are subsets of $\mathcal{C}$ and models $M, N, ..$ are elementary submodels of $\mathcal{C}$ whose cardinalities are strictly less than $\bar{\kappa}$. For simple $T$, a family of sets $\{C_i | i \in I\}$ is *$A$-independent* if for every $i \in I$ and $\bar{c}_i \in C_i$, $tp(\bar{c}_i/A \cup \bigcup \{C_j | j \neq i, j \in I\})$ does not fork over $A$. A sequence $I$ is said to be a *Morley sequence* (of $p \in S(A)$) if $I$ is $A$-indiscernible and $A$-independent, ( and $I$ is a sequence of tuples realizing $p$). Finally we note that if there is no restriction on $T$, then $T$ is arbitrary.

§**2. Lascar strong types in simple theories.** Recall that an equivalence relation $E$ (on $\mathcal{C}^n$, not necessarily definable,) is said to be *bounded* if it has strictly less than $\bar{\kappa}$ many equivalence classes. Also, for a set $A$, $E$ is said to be *$A$-invariant*, if $E(\bar{a}, \bar{b})$ implies $E(f(\bar{a}), f(\bar{b}))$ for any $A$-automorphism $f$.

**Fact 6.** ([KP]) *Suppose $E$ is an $A$-invariant bounded equivalence relation on $n$-tuples. Let $I$ be an infinite $A$-indiscernible sequence of $n$-tuples. Then $E(\bar{a}, \bar{b})$ for all $\bar{a}, \bar{b} \in I$.*

**Fact 7.** ([La2] or [KP]) $Lstp(\bar{x}/A) = Lstp(\bar{y}/A)$ *is an $A$-invariant bounded equivalence relation. (The number of classes is $\leq 2^{|T|+|A|}$.)*

**Proposition 8.** *Let $T$ be simple, and $A$ be a set. The following are equivalent.*
  (1) *For any $\bar{a}, \bar{b}$, $Lstp(\bar{a}/A) = Lstp(\bar{b}/A)$ if and only if $tp(\bar{a}/A) = tp(\bar{b}/A)$.*
  (2) *The Independence Theorem holds over $A$.*

(3) *Suppose that $tp(\bar{a}_0/A) = tp(\bar{a}_1/A)$ and $tp(\bar{a}_1/A\bar{a}_0)$ does not fork over $A$. Then there is a Morley sequence $I = \langle \bar{a}_i | i < \omega \rangle$ of $tp(\bar{a}_0/A)$.*

(4) *Suppose that $tp(\bar{a}_0/A) = tp(\bar{a}_1/A)$ and $tp(\bar{a}_1/A\bar{a}_0)$ does not fork over $A$. Then there is an $A$-indiscernible sequence $I = \langle \bar{a}_i | i < \omega \rangle$.*

*Proof.* (1)→(2) Theorem 4.

(2)→(3)→(4) Suppose that the Independence Theorem holds over $A$, and $tp(\bar{a}_1/A\bar{a}_0)$ does not fork over $A$, where $\bar{a}_1 \models tp(\bar{a}_0/A)$. We claim the following.

**Claim** *There is a sequence $I = \langle \bar{a}_i | i < \omega \rangle$ such that for each $n < \omega$,*

  (i) $tp(\bar{a}_0\bar{a}_1/A) = tp(\bar{a}_i\bar{a}_j/A)$ *for each $i < j \leq n$,*
  (ii) $\{\bar{a}_i | i \leq n\}$ *is $A$-independent.*

*Proof of Claim.* Obviously (i), (ii) hold for $n = 1$. Suppose that (i), (ii) hold for $n$ with $\langle \bar{a}_0, ... \bar{a}_n \rangle$. Let $\bar{a}_{n+1}$ be a tuple realizing a nonforking extension of $p(\bar{x}, \bar{a}_0) = tp(\bar{a}_1/A\bar{a}_0)$ over $A\bar{a}_0...\bar{a}_n$ containing $\cup_{i \leq n} p(\bar{x}, \bar{a}_i)$ (cf. Proposition 5). Then $\langle \bar{a}_0, ... \bar{a}_{n+1} \rangle$ satisfies (i), (ii) for $n + 1$. Hence the claim follows.

Now applying Ramsey's Theorem, we may assume that $I$ in the claim is $A$-indiscernible.

(4)→(1) Assume (4) holds. Suppose that $tp(\bar{a}/A) = tp(\bar{b}/A)$. Now there is $\bar{c} \models tp(\bar{a}/A)$ such that $tp(\bar{c}/A\bar{a}\bar{b})$ does not fork over $A$. Then (4) says, there are $A$-indiscernible sequences $I$, $J$ such that $\bar{c}, \bar{a} \in I$, $\bar{c}, \bar{b} \in J$. Hence by Fact 6,7, $Lstp(\bar{a}/A) = Lstp(\bar{c}/A) = Lstp(\bar{b}/A)$. □

**Corollary 9.** *Let $T$ be simple. The following are equivalent.*

(1) *The Lascar strong type over $A$ is the same as the strong type over $A$.*

(2) *The Independence Theorem holds for strong types over $A$ (cf. Theorem 4).*

(3) *Suppose that $stp(\bar{a}_0/A) = stp(\bar{a}_1/A)$ and $tp(\bar{a}_1/A\bar{a}_0)$ does not fork over $A$. Then there is an $A$-indiscernible sequence (possibly a Morley sequence) $I = \langle \bar{a}_i | i < \omega \rangle$ of $tp(\bar{a}_0/A)$.*

**Definition 10.** Let $A$ be a set, and let $\bar{a}, \bar{b}$ be tuples such that $tp(\bar{a}/A) = tp(\bar{b}/A)$. We define $d_A(\bar{a}, \bar{b})$ to be the least natural number $n(\geq 1)$ such that: there are sequences $I_1, I_2, ..., I_n$ and tuples $\bar{a} = \bar{a}_0, \bar{a}_1, ..., \bar{a}_n = \bar{b}$ such that $\bar{a}_{i-1}\frown I_i$ and $\bar{a}_i\frown I_i$ are both $A$-indiscernible for each $1 \leq i \leq n$. (If there is no such $n < \omega$, then we write $d_A(\bar{a}, \bar{b}) = \infty$.)

The following proposition says that equality of Lascar strong types is a conjunction of all bounded equivalence relations (whereas, equality of strong types is a conjunction of all definable finite equivalence relations.) A version of the proposition is already given in [KP] and [La3].

**Proposition 11.** *The following are equivalent.*

(1) $Lstp(\bar{a}/A) = Lstp(\bar{b}/A)$.

(2) $d_A(\bar{a}, \bar{b}) < \omega$.
(3) $\models E(\bar{a}, \bar{b})$ for any $A$-invariant bounded equivalence relation $E$.

*Proof.* (1)→(2) Suppose that for some model $M \supseteq A$, $tp(\bar{a}/M) = tp(\bar{b}/M)$. It suffices to show $d_A(\bar{a}, \bar{b}) = 1$. (The general case proceeds by induction.)

Let $p = tp(\bar{a}/M)$. Then there is a sequence $\langle \bar{a}_i : i < \omega \rangle$ such that $tp(\bar{a}_i/M\bar{a}\bar{b}\bar{a}_0...\bar{a}_{i-1})$ is a coheir of $p$, and is realized by $\bar{a}_{i+1}$, for $i < \omega$. Then each of $\langle \bar{a}, \bar{a}_0, \bar{a}_1, .... \rangle$ and $\langle \bar{b}, \bar{a}_0, \bar{a}_1, .... \rangle$ is an infinite $M$-indiscernible sequence. Hence $d_A(\bar{a}, \bar{b}) = 1$.

(2)→(3) Suppose $d_A(\bar{a}, \bar{b}) = 1$ (for $n$, use induction), witnessed by $I$. Then by Fact 6, for any $A$-invariant bounded equivalence relation $E$, $E(\bar{a}, \bar{c})$ and $E(\bar{c}, \bar{b})$ for some (any) $\bar{c} \in I$. Hence $E(\bar{a}, \bar{b})$.

(3)→(1) By Fact 7. □

**Example 12.** (Due to Poizat [La1]) Let us construct a model, where Lascar strong type is different from strong type. The model consists of disjoint union of $(\mathbb{R}, +, <)$ and the unit circle $C$. Also there are the additive group action of $\mathbb{R}$ on $C$ identifying $C$ with $\mathbb{R}/2\pi$; and a ternary relation $U(x, y, z)$ such that: $U(x, y, z)$ iff $x, y \in C$, $z \in \mathbb{R}$, and the length of the shorter arc from $x$ to $y$ is $< z$. We note that for each $n > 0$, there is a number $k_n$, so that whenever there are $k_n$ distinct points in $C$, then two of them should realize $U(x, y, n^{-1})$.

Now let $E$ be the equivalence relation (on $C$) defined by the conjunction of formulas $\bigwedge_{0<n\in\omega} U(x, y, n^{-1})$. By the Erdös-Rado Theorem, it is easy to see that the number of $E$-classes is $2^\omega$. But it can also be checked that for any $b, d \in C$, $stp(b/\mathbb{R}) = stp(d/\mathbb{R})$. Therefore Lascar strong type is different from strong type, (over $\mathbb{R}$).

Surprisingly, for simple $T$, equality of Lascar strong types turns out to be $\infty$-definable (i.e. type definable). Let us discuss this.

**Proposition 13.** *Let $T$ be simple. Then $Lstp(\bar{a}/A) = Lstp(\bar{b}/A)$ if and only if $d_A(\bar{a}, \bar{b}) \leq 2$.*

*Proof.* It suffices to show that $Lstp(\bar{a}/A) = Lstp(\bar{b}/A)$ implies $d_A(\bar{a}, \bar{b}) \leq 2$. (The other direction already holds, by Proposition 11.) If $Lstp(\bar{a}/A) = Lstp(\bar{b}/A)$, then Proposition 2 says there is $\bar{c}$ such that $Lstp(\bar{c}/A) = Lstp(\bar{a}/A)$ and $tp(\bar{c}/A\bar{a}\bar{b})$ does not fork over $A$. We claim that $d_A(\bar{c}, \bar{a}) = d_A(\bar{c}, \bar{b}) = 1$: By Proposition 3, there is a model $M \supseteq A$ such that $tp(\bar{c}/M) = tp(\bar{a}/M)$. Hence $d_A(\bar{c}, \bar{a}) = 1$. (See the proof of Proposition 11.(1)→(2).) For the same reason, $d_A(\bar{c}, \bar{b}) = 1$. Therefore, $d_A(\bar{a}, \bar{b}) \leq d_A(\bar{a}, \bar{c}) + d_A(\bar{c}, \bar{b}) = 2$. □

The previous proposition depends on Proposition 3. In Proposition 3, the $A$-independence of $\bar{a}, \bar{b}$ is necessary. (If $T$ is stable, then it does not need to be.) For consider the following example.

**Example 14.** Let us construct a model. The universe of the model consists of a disjoint union of countable sets $P$ and $Q$, where $P$ is a disjoint union of countable sets $U, V$ with a set isomorphism $f : U \to V$. Now there is a binary relation $R$ between $P$ and $Q$, so that $V \cup Q$ with $R$ forms the bipartite random graph, and if $x \in U, y \in Q$ then $R(x, y)$ iff $\neg R(f(x), y)$. Now $L = \{P, Q, R\}$ ($U, V, f$ are not in the language !). Then it is routine to prove that the theory $T$ of the model admits quantifier elimination for 1-formulas (but a formula $\forall z \in Q(xRz \leftrightarrow \neg yRz)$ is not quantifier eliminable), and this shows every 1-formula does not have the tree property. Hence $T$ is simple (cf. [S1, 4.2.(3)]). In $T$, it can be seen that for any set $A$, the Lascar strong type over $A$ is the same as the type over $A$ (cf. Theorem 20). Now let $a \in U$ and $b \in V$ such that $f(a) = b$. Then, as $tp(a) = tp(b)$, $Lstp(a) = Lstp(b)$, but $\{a, b\}$ is not independent over the empty set. In fact, for any model $M$, $tp(a/M) \neq tp(b/M)$ : We note that $\{R(a, M), R(b, M)\}$ forms a partition of $Q(M)$. Hence for some $c \in Q(M)$, $R(a, c)$ and $\neg R(b, c)$. This means $tp(a/M) \neq tp(b/M)$.

We further note that $d(a, b) = 1$: Let us rewrite $U = \langle a = a_0, a_1, a_2, ...\rangle$. Then $U$ and $\langle b, a_1, a_2, ...\rangle$ are both $\phi$-indiscernible.

**Theorem 15.** *Let $T$ be simple. Then for each set $A$, there is a partial type $r_A(\bar{x}, \bar{y})$ over $A$ such that $Lstp(\bar{a}/A) = Lstp(\bar{b}/A)$ if and only if $\models r_A(\bar{a}, \bar{b})$. Moreover for each $A$, $\models r_A(\bar{x}, \bar{y}) \leftrightarrow \bigwedge \{r_{\bar{c}}(\bar{x}, \bar{y}) | \bar{c} \subseteq A, \ \bar{c} \ finite\}$.*

*Proof.* By Proposition 13, we take $r_A(\bar{x}, \bar{y})$ to be a set of formulas over $A$ saying $d_A(\bar{x}, \bar{y}) \leq 2$. ($r_A(\bar{x}, \bar{y})$ says there are $\bar{x}_1, \bar{x}_2, ..., \bar{y}_1, \bar{y}_2, ...$, and $\bar{z}$ such that $\langle \bar{x}, \bar{x}_1, \bar{x}_2, ...\rangle, \langle \bar{z}, \bar{x}_1, \bar{x}_2, ...\rangle$ and $\langle \bar{z}, \bar{y}_1, \bar{y}_2, ...\rangle, \langle \bar{y}, \bar{y}_1, \bar{y}_2, ...\rangle$ are all $A$-indiscernible sequences.)

Now by compactness, $d_A(\bar{x}, \bar{y}) \leq 2$ if and only if $d_{\bar{c}}(\bar{x}, \bar{y}) \leq 2$ for all finite $\bar{c} \subseteq A$. Hence the second assertion follows. □

**Corollary 16.** *Assume $T$ to be simple. Then $Lstp(\bar{a}/A) = Lstp(\bar{b}/A)$ if and only if for all finite $\bar{c} \subseteq A$, $Lstp(\bar{a}/\bar{c}) = Lstp(\bar{b}/\bar{c})$.*

**Proposition 17.** *For simple $T$, and a set $A$, the following are equivalent.*
*(1) For any $\bar{a}, \bar{b}$, $Lstp(\bar{a}/A) = Lstp(\bar{b}/A)$ if and only if $stp(\bar{a}/A) = stp(\bar{b}/A)$.*
*(2) Assume that $E(\bar{x}, \bar{y})$ is an $\infty$-definable (i.e. type definable) bounded equivalence relation over $A$. Then $E'(\bar{x}, \bar{y}) \vdash E(\bar{x}, \bar{y})$ for some $E'(\bar{x}, \bar{y}) = \{E_i(\bar{x}, \bar{y}) | i \in I\}$, where each formula $E_i(\bar{x}, \bar{y})$ over $A$ defines a finite equivalence relation.*

*Proof.* (1)→(2) Assume (1) holds. Let us take $E'$ which defines equality of strong types over $A$, (i.e. $E'$ is a conjunction of all definable finite equivalence relations over $A$.) As the Lascar strong type over $A$ is the same as the strong type over $A$, by Proposition 11, $E'(\bar{x}, \bar{y}) \vdash E(\bar{x}, \bar{y})$.

(2)→(1) Assume (2) holds. Let $r_A(\bar{x}, \bar{y})$ define equality of Lascar strong types over $A$ (Theorem 15). Then $E'(\bar{x}, \bar{y}) \vdash r_A(\bar{x}, \bar{y})$ for some $E'(\bar{x}, \bar{y}) = \{E_i(\bar{x}, \bar{y}) | i \in I\}$, where each formula $E_i(\bar{x}, \bar{y})$ over $A$ defines a finite equivalence relation. Hence

$stp(\bar{a}/A) = stp(\bar{b}/A)$ implies $Lstp(\bar{a}/A) = Lstp(\bar{b}/A)$. On the other hand, we already know $Lstp(\bar{a}/A) = Lstp(\bar{b}/A)$ implies $stp(\bar{a}/A) = stp(\bar{b}/A)$. □

**Example 18.** (Pillay and Poizat [PPo]) Even in stable $T$, an $\infty$-definable bounded equivalence relation need not be equivalent to a conjunction of finite definable equivalence relations: The model $M$ consists of the universe $\mathbb{Q}$, and unary predicates $U_a = \{x \in \mathbb{Q} | x \leq a\}$ ($a \in \mathbb{Q}$). For each $a \in \mathbb{Q}$, let us consider the following types (over the empty set),
   (i) type $a^+$ determined by $\{U_b(x) | a < b \in \mathbb{Q}\} \cup \{\neg U_b(x) | a \not< b \in \mathbb{Q}\}$,
   (ii) type $a^-$ determined by $\{U_b(x) | a \leq b \in \mathbb{Q}\} \cup \{\neg U_b(x) | a \not\leq b \in \mathbb{Q}\}$.
Now let $E$ be the equivalence relation defined by the conjunction of the formulas $(U_a(x) \to U_b(y)) \wedge (U_a(y) \to U_b(x))$ for $a < b$. It is easy to see, for each $a \in \mathbb{Q}$, (type $a^+$)∪(type $a^-$) consists of an $E$-class.

Now if $E$ were a conjunction of (finite) definable equivalence relations $E_i$ (over $\phi$), then, for given $i$, each of (type $a^+$)∪(type $a^-$) should be contained in one of the $E_i$-classes. Quantifier elimination shows, then $E_i$ has to be trivial. ( Note that $E_i$ can not be of the form $U_a(x) \leftrightarrow U_a(y)$.) Therefore, $E$ is not a conjunction of definable equivalence relations.

§**3. Main theorems.** For the rest of this paper, $T$ will be countable.

**Proposition 19.** (Pillay, Poizat [PPo]) *Let $T$ be small (i.e. $|S(T)| \leq \omega$), and let $E(\bar{x}, \bar{y})$ be an $\infty$-definable equivalence relation over $\phi$. Suppose that $E$ is coarser than equality of n-types over $\phi$, (i.e. for any n-tuples $\bar{a}, \bar{b}$, $tp(\bar{a}) = tp(\bar{b})$ implies $E(\bar{a}, \bar{b})$). Then, there are formulas $E_i(\bar{x}, \bar{y}) \in L$ ($i \in \omega$), each of which defines an equivalence relation, such that*

$$E(\bar{x}, \bar{y}) \leftrightarrow \bigwedge_{i \in \omega} E_i(\bar{x}, \bar{y}).$$

*Proof.* Suppose that $\neg E(\bar{a}, \bar{b})$, and $p(\bar{x}) = tp(\bar{a})$, $q(\bar{y}) = tp(\bar{b})$. As $E$ is coarser than equality of types, $p(\bar{x}) \cup q(\bar{y}) \vdash \neg E(\bar{x}, \bar{y})$. Hence we can choose $\varphi(\bar{x}) \in p(\bar{x})$, $\psi(\bar{y}) \in q(\bar{y})$ such that, first, $\varphi(\bar{x}) \wedge \psi(\bar{y}) \vdash \neg E(\bar{x}, \bar{y})$ and secondly, the set $X$, defined by $\neg(\varphi(\bar{x}) \vee \psi(\bar{x}))$, has minimum Cantor-Bendixson rank and degree. We claim $X = \phi$:

If not, then $CB(X) = (\alpha, n)$. Now there is $p'(\bar{y}) \in S(\phi)$ containing $\neg(\varphi(\bar{y}) \vee \psi(\bar{y}))$ such that the $CB$-rank of $p'(\bar{y})$ is $\alpha$. Let $\bar{c} \models p'(\bar{y})$. Then either $\neg E(\bar{a}, \bar{c})$ or $\neg E(\bar{b}, \bar{c})$. Without loss of generality, we assume $\neg E(\bar{a}, \bar{c})$. Then in exactly the same manner, we obtain $\sigma(\bar{y}) \in p'(\bar{y})$ such that $\varphi(\bar{x}) \wedge \sigma(\bar{y}) \vdash \neg E(\bar{x}, \bar{y})$. Now it can be seen that $CB(\neg(\varphi(\bar{x}) \vee (\sigma(\bar{y}) \vee \psi(\bar{y})))) < (\alpha, n)$. This contradicts the choice of $\varphi$ and $\psi$. Hence the claim follows.

Now we note that $\{\varphi(\bar{x}), \psi(\bar{x})\}$ forms a partition of a (saturated) model $M$. Let us rewrite $\varphi(\bar{x})$ as $\varphi_{\bar{a}\bar{b}}(\bar{x})$. It follows that $E(\bar{x}, \bar{y})$ is a conjunction of formulas $\varphi_{\bar{a}\bar{b}}(\bar{x}) \leftrightarrow \varphi_{\bar{a}\bar{b}}(\bar{y})$, for $\bar{a}\bar{b} \in M$ with $\neg E(\bar{a}, \bar{b})$. □

**Theorem 20.** *Let $T$ be simple, and small. Then Lascar strong type is the same as strong type, over any set $A$.*

*Proof.* We note that $stp(\bar{a}/A) = stp(\bar{b}/A)$ if and only if for all finite $\bar{c} \subseteq A$, $stp(\bar{a}/\bar{c}) = stp(\bar{b}/\bar{c})$. Hence, by Corollary 16, it suffices to prove the theorem for finite $A$. Moreover, if $A$ is finite, then as $T_A$ is again small, we may assume $A = \phi$. Now, there is a type $r(\bar{x}, \bar{y})$, which defines Lascar strong types (over $\phi$.) Suppose that $stp(\bar{c}) = stp(\bar{d})$. We want to show $Lstp(\bar{c}) = Lstp(\bar{d})$.

**Claim** *Let $p(\bar{x}) = tp(\bar{c})$. Then there are formulas $R_i(\bar{x}, \bar{y}) \in L$ $(i \in \omega)$, each of which defines an equivalence relation, such that*

$$p(\bar{x}) \cup p(\bar{y}) \cup r(\bar{x}, \bar{y}) \leftrightarrow p(\bar{x}) \cup p(\bar{y}) \cup \{R_i(\bar{x}, \bar{y}) | i \in \omega\} \quad (1).$$

*Proof of Claim* For the proof of this claim, we borrow the technique Lascar used to prove his result on "$\wedge$-basic subgroup" ([La3]).

Now let $E(\bar{x}, \bar{y})$ be an equivalence relation over $\bar{c}$ defined as $E(\bar{x}, \bar{y})$ iff $\exists \bar{z}(tp(\bar{c}\bar{x}) = tp(\bar{z}\bar{y}) \wedge r(\bar{c}, \bar{z}))$. Then if $tp(\bar{a}/\bar{c}) = tp(\bar{b}/\bar{c})$, then obviously $E(\bar{a}, \bar{b})$. Hence by the preceding proposition,

$$E(\bar{x}, \bar{y}) \leftrightarrow \bigwedge_{i \in \omega} E_i(\bar{x}, \bar{y}; \bar{c}) \quad (*)$$

for some formulas $E_i(\bar{x}, \bar{y}; \bar{c})$ over $\bar{c}$, each of which defines an equivalence relation. Now, for each $i \in \omega$, let

$$R_i(\bar{x}, \bar{y}) \equiv \forall \bar{u}(E_i(\bar{u}, \bar{x}, \bar{x}) \leftrightarrow E_i(\bar{u}, \bar{y}, \bar{y})).$$

We claim that formulas $R_i(\bar{x}, \bar{y})$ are the desired formulas: It suffices to show $p(\bar{y}) \models r(\bar{c}, \bar{y}) \leftrightarrow \wedge_{i \in \omega} R_i(\bar{c}, \bar{y})$. Suppose that $\bar{a} \models p(\bar{x})$. Now let $r(\bar{c}, \bar{a})$. Then, by (*) and the definition of $E(\bar{x}, \bar{y})$, if $tp(\bar{a}\bar{u}) = tp(\bar{c}\bar{u}')$, then $E_i(\bar{u}, \bar{u}'; \bar{c})$ for each $i$. Hence $E_i(\bar{u}, \bar{c}, \bar{c}) \leftrightarrow E_i(\bar{u}', \bar{c}, \bar{c}) \leftrightarrow E_i(\bar{u}, \bar{a}, \bar{a})$. Therefore $R_i(\bar{c}, \bar{a})$ for each $i \in \omega$.

Conversely, if $R_i(\bar{c}, \bar{a})$ for each $i$, then as $E_i(\bar{a}, \bar{a}, \bar{a})$, $E_i(\bar{a}, \bar{c}, \bar{c})$ for each $i$. Thus $E(\bar{a}, \bar{c})$ and $r(\bar{a}, \bar{c})$. $\square$

Let us continue to prove the theorem. Now, for each $i \in \omega$ and each formula $\psi(\bar{x})$ in $p(\bar{x})$,

$$(\psi(\bar{x}) \wedge \psi(\bar{y}) \wedge R_i(\bar{x}, \bar{y})) \vee (\neg \psi(\bar{x}) \wedge \neg \psi(\bar{y})) \stackrel{let}{=} R_i^\psi(\bar{x}, \bar{y})$$

defines an equivalence relation. We note that, each equivalence class of $R_i^\psi$ is either $\neg \psi(\bar{x})$, or $\psi(\bar{x}) \wedge$ (an equivalence class of $R_i$). Moreover, for each $i \in \omega$, there is some formula $\varphi_i(\bar{x})$ in $p(\bar{x})$ such that $R_i^{\varphi_i}(\bar{x}, \bar{y})$ defines a finite equivalence relation. ( If not, then by (1), $r(\bar{x}, \bar{y})$ can not be bounded; a contradiction.) Now,

$$p(\bar{x}) \cup p(\bar{y}) \cup \{R_i^{\varphi_i}(\bar{x}, \bar{y}) | i \in \omega\} \rightarrow \{\varphi_i(\bar{x}) \wedge \varphi_i(\bar{y}) \wedge R_i^{\varphi_i}(\bar{x}, \bar{y}) | i \in \omega\}$$
$$\rightarrow \{\varphi_i(\bar{x}) \wedge \varphi_i(\bar{y}) \wedge R_i(\bar{x}, \bar{y}) | i \in \omega\}. \quad (2)$$

A NOTE ON LASCAR STRONG TYPES IN SIMPLE THEORIES 9Finally, since $stp(\bar{c}) = stp(\bar{d})$, $p(\bar{c}) \cup p(\bar{d}) \cup \{R_i^{\varphi_i}(\bar{c},\bar{d})|i \in \omega\}$. Hence, by (2), $\{R_i(\bar{c},\bar{d})|i \in \omega\}$, and by (1), $r(\bar{c},\bar{d})$. Therefore $Lstp(\bar{c}) = Lstp(\bar{d})$. □

Theorem 20 recovers the fact that the Independent Theorem holds over an algebraically closed set (in $\mathcal{C}^{eq}$) in any smoothly approximable structures ([CH], [KaLM]).

**Corollary 21.** *If $T$ is simple and $\omega$-categorical, then Lascar strong type is the same as strong type, over any set.*

For a final remark, we define terminology from [PPo]. For each $\infty$-definable equivalence relation $E(\bar{x},\bar{y})$ over $A$, we write $E^*(\bar{x},\bar{y})$ for a corresponding equivalence relation such that,
$$\models E^*(\bar{a},\bar{b}) \text{ iff there is } \bar{b}' \text{ such that } \bar{b}' \models tp(\bar{b}/A) \text{ and } \models E(\bar{a},\bar{b}').$$
It is easy to check that, $E^*$ is again an $\infty$-definable equivalence relation over $A$, and is coarser than equality of types over $A$.

Pillay and Poizat ([PPo]) proved that for any $T$, an $\infty$-definable equivalence relation $E(\bar{x},\bar{y})$ over $A$ is equivalent to a conjunction of definable equivalence relations over $A$ if and only if, so are $E^*(\bar{x},\bar{y})$ and each restriction of $E(\bar{x},\bar{y})$ to each $n$-type $p \in S(A)$. Now if $T$ is small, then for given $E$ (over $\phi$), $E^*$ is already equivalent to a conjunction of definable equivalence relations (Theorem 19). Moreover, by exactly the same argument in the proof of Claim in Theorem 20, so is each restriction of $E$ to each $p \in S(\phi)$. Hence we obtain the following theorem. This theorem also proves Theorem 20.

**Theorem 22.** *Let $T$ be small, and let $A$ be a finite set. Then any $\infty$-definable equivalence relation over $A$ is equivalent to a conjunction of definable equivalence relations over $A$.*

## References

[CH] G. Cherlin and E. Hrushovski, 'Lie coordinatised structures', in preparation.
[KaLM] W. M. Kantor, M. W. Liebeck, and H. D. Macpherson, '$\aleph_0$-categorical structures smoothly approximated by finite structures', *Proc. London Math. Soc.* (3) 59 (1989) 439-463.
[K1] B. Kim, 'Forking in simple unstable theories', to appear in *J. of London Math. Soc.*
[K2] B. Kim, 'Recent results on simple first order theories', to appear in *the Blaubeuren proceedings* London Math. Soc. Lecture Note Series (Cambridge University Press).
[KP] B. Kim and A. Pillay, 'Simple theories', to appear in *Ann. of Pure and Applied Logic.*
[La1] D. Lascar, *The group of automorphisms of a relational saturated structure*, ( *Proceedings, Banff, 1991*) NATO ASI Series C: Math. Phys. Sci. 411 (Kluwer Acad. Publ., Dordrecht, 1993) 225-236.
[La2] D. Lascar, 'On the category of models of a complete theory', *J. of Symbolic Logic* 47 (1982) 249-266.
[La3] D. Lascar, 'Recovering the action of an automorphism group', to appear in *Proceedings, Keele.*
[PPo] A. Pillay and B. Poizat, 'Pas d'imaginaires dans l'infini!' *J. of Symbolic Logic* 52 (1987) 400-403.

The Fields Institute, 222 College Street, Toronto, Ontario, Canada, M5T 3J1
*E-mail address*: bkim@fields.utoronto.ca